\documentclass{amsart}
\usepackage{graphicx}
\usepackage{eucal}
\usepackage{amssymb}
\usepackage{mathrsfs}
\usepackage{amssymb}
\usepackage{amsmath}
\usepackage{picture}
\usepackage{graphics}
\usepackage{mathrsfs}
\usepackage{color,cite,graphicx}
\usepackage{tikz}
\definecolor{gray}{gray}{0.4}
\usepackage[all]{xy}
\usepackage{enumerate}
 \usepackage{rotating}
 \usepackage{caption}
\title[Order eight]{Symmetries of order eight on K3 surfaces without high genus curves in the fixed locus}
\author{Dima Al Tabbaa}
\address{Universit\'e de Poitiers
Laboratoire de Math\'ematiques et Applications
UMR 7348 du CNRS
B\^at. H3 - Site du Futuroscope
TSA 61125
11 bd Marie et Pierre Curie
86073 POITIERS Cedex 9
France}
\email{Dima.Al.Tabbaa@math.univ-poitiers.fr}

\author{Annalisa Grossi}
\address{Dipartimento di Matematica, Universit\`{a} di Bologna, Piazza di Porta S. Donato, 5, 40126 Bologna, Italy and Universit\"at Augsburg, Institut f\"ur Mathematik,  Universit\"atsstr. 14, 86159 Augsburg, Germany}
\email{annalisa.grossi3@unibo.it, 	annalisa.grossi@math.uni-augsburg.de}
\urladdr{https://www.unibo.it/sitoweb/annalisa.grossi3}
\author{Alessandra Sarti}
\address{Universit\'e de Poitiers
Laboratoire de Math\'ematiques et Applications
UMR 7348 du CNRS
B\^at. H3 - Site du Futuroscope
TSA 61125
11 bd Marie et Pierre Curie
86073 POITIERS Cedex 9
France}
\email{sarti@math.univ-poitiers.fr}
\urladdr{http://www-math.sp2mi.univ-poitiers.fr/~sarti/}

\newtheorem{theorem}{Theorem}[section]
\newtheorem{pro}[theorem]{Proposition}
\newtheorem{lemma}[theorem]{Lemma}

\newtheorem{definition}[theorem]{Definition}

\newtheorem{remark}[theorem]{Remark}
\DeclareMathOperator{\Pic}{Pic}

\DeclareMathOperator{\Fix}{Fix}

\DeclareMathOperator{\rk}{rk}
\DeclareMathOperator{\Aut}{Aut}

\DeclareMathOperator{\NS}{NS}

\newcommand{\IP}{\mathbb{P}}

\usepackage{multicol}

\keywords{non--symplectic automorphisms, K3 surfaces}
\thanks{}

\begin{document}

\maketitle

\begin{abstract}
In this paper we classify non--symplectic automorphisms of order eight on complex K3 surfaces in case that the fourth power of the automorphism has only rational curves in its fixed locus. We show that the fixed locus is the disjoint union of a rational curve and ten isolated points or it consists in four isolated fixed points. We give examples corresponding to the case with a rational curve in the fixed locus and to the case with only isolated points in the fixed locus. 
\end{abstract}

\section*{Introduction}

In this paper we investigate purely non--symplectic automorphisms of order eight on a complex $K3$ surface under certain assumptions on the fixed locus. These automorphisms act non--trivially on $H^{2,0}(X) \simeq\mathbb{C} \cdot \omega_{X}$ i.e. they multiply the non--degenerate holomorphic 2--form by a primitive 8th root of unity. 
The study of non--symplectic automorphisms of prime order was completed by several authors: Nikulin  \cite{nikulinfactor}, Artebani, Sarti and Taki\cite{AS3}, \cite{ast} and \cite{takiauto}. If the automorphism is not of prime order the setting is more complicated. Indeed, in this situation 
the purely non--symplectic automorphism of order eight does not admit a trivial action on the N\'eron--Severi group  of the generic K3 surface as it does in the case of the prime order \cite[Section 11]{dolgachev2007moduli}. In the paper \cite{Taki}, Taki studies the case when the order of the automorphism is a prime power and the action is trivial on the N\'eron--Severi group. If we consider non--symplectic, non--trivial automorphisms of order $2^b$, then by results of Nikulin we have $1\leq b\leq 5$. Further results can be found in a paper by Sch\"{u}tt \cite{matthias} in the case of automorphisms of a two--power order and in a paper by Artebani and Sarti \cite{ast} in the case of order four. Recently in \cite{tabbaa2014classification} Al Tabbaa, Sarti and Taki completed the study for purely non--symplectic automorphisms of order $16$. In \cite{tabbaa2016order} Al Tabbaa and Sarti studied the case of order eight automorphisms under the assumption that their fourth power $\sigma^{4}$ acts as the identity on the N\'eron--Severi group and Fix($\sigma^{4}$) contains an elliptic curve. \\
This paper deals with purely non--symplectic automorphisms of order eight on K3 surfaces under the assumption that their fourth power $\sigma^{4}$ is the identity on the N\'eron--Severi group. This corresponds to the situation for the generic K3 surface in the moduli space of K3 surfaces with a purely non--symplectic automorphism of order eight  and fixed action on the second cohomology group with integer coefficients, see \cite[Section 10]{dolgachev2007moduli}. The fixed locus Fix($\sigma$) of such an automorphism $\sigma$ is the disjoint union of smooth curves and points. We give a complete classification of the fixed locus and of some lattice invariants related to $\sigma$ in the case that Fix($\sigma^{4}$) is not empty and contains only rational curves. Let $X$ be a K3 surface, $\omega_{X}$ a generator of $H^{2,0}(X)$ and $\sigma \in \Aut(X)$ such that $\sigma^{\ast}(\omega_{X})=\zeta_{8}\omega_{X}$, 
where $\zeta_8$ denotes a primitive 8th root of unity. We denote by $k_{\sigma}$ the number of smooth rational curves fixed by $\sigma$ and by $N_{\sigma}$ the numbers of isolated points in Fix($\sigma$). We denote furthermore by $S(\sigma^{4})$ the invariant lattice for the action of $\sigma^4$ on the second cohomology group with integer coefficients. We prove the following result:

\begin{theorem}\label{main}
	Let $\sigma$ be a purely non--symplectic automorphism of order 8 on a K3 surface $X$ with $Pic(X)=S(\sigma^{4})$. Suppose that $\Fix(\sigma^{4})$ is not empty and it is the union of smooth rational curves. Then $k_{\sigma} \in \{0, 1\}$ moreover
	\begin{itemize}
		\item[$\bullet$]if $k_{\sigma}=1$ then $N_{\sigma}=10$,
		\item[$\bullet$]if $k_{\sigma}=0$ then $N_{\sigma}=4$.
	\end{itemize}
\end{theorem}

The paper is organized as follows. In Section \ref{basic_s} we recall basic facts on purely non--symplectic automorphisms acting on K3 surfaces. 
In Section \ref{fixed_s} we give several properties of the fixed locus of such automorphisms in particular in the case that the automorphism has even order. The Section \ref{classification_s} is devoted to the proof of our main theorem \ref{main}. In Table \ref{tabella thm} we list in detail all the possibilities for the fixed locus and for some invariants characterizing the action 
of the automorphism on the second cohomology group with integer coefficients. We have in total four possibilities:  one for $k_{\sigma}=1$ and three for $k_{\sigma}=0$. By using elliptic fibrations we give examples for $k_{\sigma}=1$ and for $k_{\sigma}=0$ in the case 
that the number of fixed curves by $\sigma^2$ is three, this is done in Section \ref{example_s} after recalling in Section \ref{elliptic_s} some basic notions on elliptic K3 surfaces. We do not know if the remaining two cases in Table \ref{tabella thm} exist. In any case by using the invariants of the Table \ref{tabella thm} one computes that the rank of the N\'eron--Severi group is $18$ or $14$ so that by e.g. \cite[Lemma 1.5]{machida_oguiso} all these K3 surfaces admits an elliptic fibration. This could certainly be helpful in finding the remaining examples.\\

\textbf{Acknowledgements}: The second author was supported by the Laboratoire International LIA LYSM to visit the University of Poitiers
in May 2019. We thank Paola Comparin for useful comments and the anonymous referee for a very careful reading of the paper and for the several useful remarks, which improved the paper.

\section{Basic facts}\label{basic_s}
Let $X$ be a K3 surface and $\sigma \in \Aut(X)$ a non--symplectic automorphism of order $8$. We assume that $\sigma^{\ast}(\omega_{X})= \zeta_8 \omega_{X}$ where $\zeta_8$ is a primitive 8th root of unity. Such a $\sigma$ is called \emph{purely non--symplectic}, for simplicity we just call it $\emph{non--symplectic}$, always meaning that the action on the holomorphic two--form is by a primitive 8th root of unity. \\
We denote by $r_{\sigma^{j}}, l_{\sigma^{j}}, m_{\sigma^{j}}$ and $m_1$ for $j=1,2,4$ the rank of the eigenspace of $(\sigma^{j})^{\ast}$ in $H^{2}(X, \mathbb{C})$ relative to the eigenvalues $1$, $-1$, $i$ and $\zeta_{8}$ respectively (clearly $m_{\sigma^{4}}=0)$. We recall the invariant lattice:
$$ S(\sigma^{j})=\{x \in H^{2}(X, \mathbb{Z}) | (\sigma^{j})^{\ast}(x)=x \},$$
and its orthogonal complement 
$$ T(\sigma^{j})=S(\sigma^{j})^{\perp} \cap H^{2}(X, \mathbb{Z}).$$
Since the automorphism acts purely non--symplectically, $X$ is projective, see \cite[Theorem 3.1]{Nikulin1}, so that if we denote $\rk S(\sigma^{j})=r_{\sigma^{j}}$, we have that $r_{\sigma^{j}} > 0$ for all $j=1,2,4$, in fact one can always find an invariant ample class. On the other hand, one can easily show that $S(\sigma^{j}) \subset \Pic(X)$ for $j=1,2,4$ so that the transcendental lattice satisfies $T_{X} \subseteq T(\sigma^{j})$ for $j=1,2,4$. 

\begin{remark}\label{rmk 0}
	It is a straightforward computation that the invariants $r_{\sigma^{j}}, l_{\sigma^{j}}, m_{\sigma^{j}}$ and $m_1$ with $j=1,2,4$ satisfy the following relations: 
	\begin{center}
		\begin{multicols}{2}
			$r_{\sigma^{2}}=r_{\sigma} + l_{\sigma} $; \\ $l_{\sigma^{2}}=2m_{\sigma}$; \\ $m_{\sigma^{2}}=2m_1$. \\
			$r_{\sigma^{4}}=r_{\sigma}+l_{\sigma}+2m_{\sigma}$; \\ $l_{\sigma^{4}}=4m_1$;
		\end{multicols}
	\end{center}
	We remark that the invariants $l_{\sigma^{2}}$ and $m_{\sigma^{2}}$ are even numbers.
\end{remark}  

The moduli space for $K3$ surfaces carrying a non--symplectic automorphism of even order $n>2$, with a given action on the $K3$ lattice is known to be a complex ball quotient of dimension $q-1$ where $q$ is the dimension of the eigenspace $V$ of $\sigma^{*}$ in $H^{2}(X, \mathbb{C})$ relative to the eigenvalues $\zeta_{n}=e^{\frac{2\pi i }{n}}$, see \cite[\S 11]{dolgachev2007moduli}. The complex ball is given by:
$$ B=\{[w] \in \mathbb{P}(V) : (w, \overline{w}) >0 \}.$$
If $n$ is even $V$ is the $\zeta_{n}$ eigenspace of $\sigma^{*}$ in $T(\sigma^{n/2}) \otimes \mathbb{C}$. This implies that the N\'eron--Severi group of a $K3$ surface corresponding to the generic point in the moduli space equals $S(\sigma^{n/2})$ see \cite[Theorem 11.2]{dolgachev2007moduli}.

\section{The fixed locus}\label{fixed_s}

We denote by Fix($\sigma^{j}$), $j=1,2,4$ the fixed locus of the automorphism $\sigma^{j}$ :
\begin{center}
	Fix($\sigma^{j})=\{x \in X | \  \sigma^{j}(x)=x \}. $
\end{center}
Clearly Fix($\sigma$)$ \subset$ Fix($\sigma^{2}$) $\subset$ Fix($\sigma^{4}$). To describe the fixed locus of order $8$ non--symplectic automorphisms we start recalling the following result about non--symplectic involutions, see \cite[Theorem 4.2.2]{nikulinfactor}.

\begin{theorem}\label{Nik invo}
	Let $\tau$ be a non--symplectic involution on a K3 surface $X$. The fixed locus of $\tau$ is either empty, the disjoint union of two elliptic curves or the disjoint union of a smooth curve of genus $g \geq 0$ and $k$ smooth rational curves.
	Moreover, its fixed lattice $S(\tau) \subset \Pic(X)$ is a $2$--elementary lattice with determinant $2^a$ such that: 
	\begin{itemize}
		\item [$\bullet$] $S(\tau) \cong U(2) \oplus E_8(2)$ iff the fixed locus of $\tau$ is empty;
		\item[$\bullet$] $S(\tau) \cong U \oplus E_8(2)$ iff $\tau$ fixes two elliptic curves; 
		\item[$\bullet$] $2g=22-rkS(\tau) - a$ and $2k=rkS(\tau) -a$ otherwise.
	\end{itemize}
\end{theorem}
Since $S(\tau)$ is $2$--elementary its discriminant group $A_{S(\tau)} = S(\tau)^{\vee}/S(\tau) \simeq (\mathbb{Z}/2\mathbb{Z})^{\oplus a}$, $a \in \mathbb{Z}_{>0}$. We introduce the invariant $\delta_{S(\tau)}$ of $S(\tau)$ by putting $\delta_{S(\tau)}=0$ if $x^{2} \in \mathbb{Z}$ for any $x \in A_{S(\tau)}$ and $\delta_{S(\tau)}=1$ otherwise. By \cite[Theorem 3.6.2]{nikulinintegral}, and \cite[\S 1]{rudakov1981surfaces} $S(\tau)$ is uniquely determined by the invariant $\delta_{S(\tau)}$, by the rank, the signature and the invariant $a$. The situation is summarized in Figure \ref{ord2} from \cite[\S 4]{nikulindiscrete}.
\begin{figure}[h]
	$$\begin{array}{cccccccccccccccccccccccr}
	&\ &\ &&&& &&\ &\ &&&&&&&&&&&&& \bullet\ \ \delta_{S(\tau)}=1\\
	&\ &\ &&&& &&\ &\ &&&&&&&&&&&&&\ast\ \ \delta_{S(\tau)}=0
	\end{array}$$
	\vspace{-1.1cm}
	
	\begin{tikzpicture}[scale=.43]
	\filldraw [black]
	(1,1) circle (1.5pt)  node[below=-0.55cm]{10}
	(2,0) node[below=-0.20cm]{*}
	(2,2) circle (1.5pt) node[below=-0.5cm]{9} node[below=-0.15cm]{*}
	(3,1) circle (1.5pt)
	(3,3) circle (1.5pt)node[below=-0.5cm]{8}
	(4,2) circle (1.5pt)
	(4,4) circle (1.5pt)node[below=-0.5cm]{7}
	(5,3) circle (1.5pt)
	(5,5) circle (1.5pt)node[below=-0.5cm]{6}
	(6,4) circle (1.5pt)node[below=-0.15cm]{*}
	(6,2)    node[below=-0.23cm]{*}
	(6,6) circle (1.5pt)node[below=-0.5cm]{5}
	(7,3) circle (1.5pt)
	(7,5) circle (1.5pt)
	(7,7) circle (1.5pt)node[below=-0.5cm]{4}
	(8,2) circle (1.5pt)
	(8,4) circle (1.5pt)
	(8,6) circle (1.5pt)
	(8,8) circle (1.5pt)node[below=-0.5cm]{3}
	(9,1) circle (1.5pt)
	(9,3) circle (1.5pt)
	(9,5) circle (1.5pt)
	(9,7) circle (1.5pt)
	(9,9) circle (1.5pt)node[below=-0.5cm]{2}
	(10,0)  node[below=-0.20cm]{*}
	(10,2) circle (1.5pt)node[below=-0.15cm]{*}
	(10,4) circle (1.5pt)node[below=-0.15cm]{*}
	(10,6) circle (1.5pt)node[below=-0.15cm]{*}
	(10,8) circle (1.5pt)node[below=-0.15cm]{*}
	(10,10) circle (1.5pt) node[below=-0.5cm]{1}node[below=-0.15cm]{*}
	(11,1) circle (1.5pt)
	(11,3) circle (1.5pt)
	(11,5) circle (1.5pt)
	(11,7) circle (1.5pt)
	(11,9) circle (1.5pt)
	(11,11) circle (1.5pt) node[below=-0.5cm]{0}
	(12,2) circle (1.5pt)
	(12,4) circle (1.5pt)
	(12,6) circle (1.5pt)
	(12,8) circle (1.5pt)
	(12,10) circle (1.5pt) node[below=-0.5cm]{1}
	(13,3) circle (1.5pt)
	(13,5) circle (1.5pt)
	(13,7) circle (1.5pt)
	(13,9) circle (1.5pt) node[below=-0.5cm]{2}
	(14,2)  node[below=-0.15cm]{*}
	(14,4) circle (1.5pt)node[below=-0.15cm]{*}
	(14,6) circle (1.5pt)node[below=-0.15cm]{*}
	(14,8) circle (1.5pt) node[below=-0.5cm]{3}
	(15,3) circle (1.5pt)
	(15,5) circle (1.5pt)
	(15,7) circle (1.5pt) node[below=-0.5cm]{4}
	(16,2) circle (1.5pt)
	(16,4) circle (1.5pt)
	(16,6) circle (1.5pt) node[below=-0.5cm]{5}
	(17,1) circle (1.5pt)
	(17,3) circle (1.5pt)
	(17,5) circle (1.5pt)node[below=-0.5cm]{6}
	(18,0)  node[below=-0.2cm]{*}
	(18,2) circle (1.5pt)node[below=-0.15cm]{*}
	(18,4) circle (1.5pt) node[below=-0.5cm]{7}node[below=-0.15cm]{*}
	(19,1) circle (1.5pt)
	(19,3) circle (1.5pt) node[below=-0.5cm]{8}
	(20,2) circle (1.5pt) node[below=-0.5cm]{9}
	;
	\draw plot[mark=*] file {data/ScatterPlotExampleData.data};
	\draw[->] (0,0) -- coordinate (x axis mid) (22,0);
	\draw[->] (0,0) -- coordinate (y axis mid)(0,12);
	\foreach \x in {0,1,2,3,4,5,6,7,8,9,10,11,12,13,14,15,16,17,18,19,20}
	\draw [xshift=0cm](\x cm,0pt) -- (\x cm,-3pt)
	node[anchor=north] {$\x$};
	\foreach \y in {1,2,3,4,5,6,7,8,9,10,11}
	\draw (1pt,\y cm) -- (-3pt,\y cm) node[anchor=east] {$\y$};
	\node[below=0.2cm, right=4.5cm] at (x axis mid) {$\rk(S(\tau))$};
	\node[left=0.5cm, below=-2.7cm, rotate=90] at (y axis mid) {$a$};
	\draw[<-, blue](0.1,0.1)-- node[below=2cm,left=2cm]{$g$} (11,11);
	\draw[<-, red](21.5,0.5)-- node[below=2cm,right=2.1cm]{$k$} (11,11);
	\end{tikzpicture}
	\caption{Order 2}\label{ord2}
\end{figure}

We recall a result about non--symplectic automorphisms of order four on a K3 surface. These results are discussed in \cite{ASorder4}, see also
the Appendix of \cite{tabbaa2014classification}. 

\begin{theorem}\label{rational}
	Let $X$ be a K3 surface and $\gamma$ be a purely non--symplectic automorphism of order four acting on it, with $\Pic(X)=S(\gamma^2)$.
	If $\Fix(\gamma)$ contains a smooth rational curve and all curves fixed by $\gamma^2$ are rational,
	then the invariants associated to $\gamma$ are as in Table \ref{g=0}. All cases in the table do exist. Here $m$ denotes the multiplicity of the eigenvalue $i$; $r$ denotes the multiplicity of the eigenvalue $1$; $l$ the multiplicity of the eigenvalue $-1$ and $a_\gamma$ the number of rational curves exchanged by $\gamma$ and fixed by $\gamma^2$. Moreover $N$ denotes the number of isolated fixed points and $k$ denotes the number of fixed rational curves.
\end{theorem}

\begin{table}[ht]
	$$
	\begin{array}{ccc|ccccccl}
	m & r & l &  N & k & a_\gamma  \\
	\hline
	4 & 10 & 4 &  6 & 1 & 0   \\
	\hline
	3 & 13 & 3 &  8 &2 & 0 \\
	& 11 & 5 &  6& 1 & 1 \\
	\hline
	2 & 16 & 2 &  10 &3 & 0 \\
	& 14 & 4 &  8 &2 & 1 \\
	& 12 & 6 &  6& 1 & 2 \\
	\hline
	1 & 19 & 1 & 12 & 4 & 0  \\
	& 13 & 7 &  6 & 1 & 3  \\
	\end{array}
	$$
	\vspace{0.2cm}
	
	\caption{The case $g=0$}\label{g=0}
\end{table}

\begin{remark}\label{rmk 1}
For each $p\in$ $\Fix(\sigma)$, there exists a $\sigma^{4}$--fixed smooth rational curve $R$, such that 
$p\in R$. 
\end{remark}
With the notation and the assumptions of the previous remark we have
\begin{lemma}\label{rmk 2} The curve $R$ is $\sigma$--invariant.
\end{lemma}
\proof 
First of all we notice that since $R$ is fixed by $\sigma^{4}$ then also $\sigma(R)$ is fixed by $\sigma^{4}$. If $R$ was not $\sigma$--invariant this means that $R$ is sent to another rational curve $\sigma(R)$ which meets $R$ in $p$.  This is impossible since these two curves are fixed by $\sigma^{4}$ and the fixed locus of an involution is smooth by Theorem $\ref{Nik invo}$.
\endproof

Observe now that in a neighborhood of a fixed point $p$ of $\sigma$ the action can be linearized, see e.g. \cite[Section 5]{Nikulin1}. We can assume $p=(0,0)$ and we can find $z_1$ and $z_2$ local coordinates in a neighborhood of $p$ such that the everywhere non--degenerate holomorphic 2--form can be written as $dz_{1} \wedge dz_2$. Since the automorphism is of finite order locally it can be diagonalized and we know that $\sigma^{*}(dz_1 \wedge dz_2)=\zeta_8(dz_1 \wedge dz_2)$, for this reason the product of the eigenvalues with respect to $z_1$ and $z_2$ is equal to $\zeta_8$. We get the following possibilities up to permutation of the coordinates (but this does not play any role in the classification):
\begin{center}
\begin{multicols}{4}
	$A_{1,0}= 
	\begin{pmatrix}
	\zeta_8 & 0 \\
	0     &      1 
	\end{pmatrix}$,	
	
		$A_{2,7}=
	\begin{pmatrix}
	i & 0 \\
	0     &  \zeta_8^{7}      
	\end{pmatrix}$,

	$A_{3,6}=
	\begin{pmatrix}
	\zeta_8^{3} & 0 \\
	0     &      \zeta_8^6
	\end{pmatrix}$,	
	
	$A_{4,5}=
	\begin{pmatrix}
	-1 & 0 \\
	0     &      \zeta_8^5 
	\end{pmatrix}$.
\end{multicols}
\end{center}

In the first case the point belongs to a smooth fixed curve, since we have an eigenvalue which is equal to $1$. In the other three cases it is an isolated fixed point. We say that an isolated point $p \in $ $\Fix(\sigma)$ is of type $(t,s)$ if the local action at $p$ is given by $A_{t,s}$. We denote by $n_{t,s}$ the number of isolated fixed points by $\sigma$  with matrix $A_{t,s}$. \\

We further denote by $N_{\sigma^{j}}$, $k_{\sigma^{j}}$, $j=1,2,4$ the number of isolated points and smooth rational curves in $\Fix(\sigma^{j})$. We observe that $N_{\sigma^{4}}=0$ since $\sigma^{4}$ only fixes curves or it is empty as explained in Theorem \ref{Nik invo}.  We recall \cite[Proposition 2.2] {tabbaa2016order} :

\begin{pro}\label{order 8}
	Let $\sigma$ be a non--symplectic automorphism of order 8 acting on a K3 surface $X$. Then $\Fix(\sigma)$ is never empty nor it can be the union of two smooth elliptic curves. It is the disjoint union of smooth curves and $N_{\sigma} \geq 2$ isolated points. We denote $\alpha=\sum\limits_{K \subset \Fix(\sigma)} (1-g(K))$ and the following relations hold:
	\begin{center}
		$n_{2,7} +n_{3,6}=2+4 \alpha,\,\, n_{4,5} + n_{2,7}-n_{3,6}=2+2\alpha,\,\, N_{\sigma}=2+r_{\sigma}-l_{\sigma}-2\alpha.$ \\
	\end{center}

\end{pro}
\vspace{0.5cm}

The fixed locus of such an automorphism $\sigma$ is then 
\begin{equation}\label{luogofisso}
\Fix(\sigma)=C \cup R_1 \cup \dots \cup R_k \cup \{p_1, \dots, p_N\} 
\end{equation}
where $C$ is a smooth curve of genus $g \geq 0$ and $R_j$, $j=1,\ldots,k$ are disjoint smooth rational curves. \\

We recall the following remarks and lemmas which are important in the study of the fixed locus of $\sigma$.
\begin{remark}
A non--symplectic automorphism $\sigma$ of order $8$ acts on a set of smooth rational curves of $X$ which are fixed by $\sigma^{4}$ either trivially, i.e. each smooth rational curve is $\sigma$--invariant, possibly pointwise fixed by $\sigma$, or it exchanges smooth rational curves two by two, or  finally $\sigma$ permutes four rational curves. 
\end{remark}
\begin{lemma}
Four cyclically permuted smooth rational curves by a non--symplectic automorphism $\sigma$ of order $8$ on a K3 surface $X$ are either $\sigma^{4}$--invariant (not pointwise fixed), or pointwise fixed by $\sigma^{4}$.
\end{lemma}

\proof 
We can prove it simply as follows. Let $C_{i}$, $i \in\{ 1, \ldots, 4\}$ be four smooth rational curves such that $\sigma(C_{i}) = C_{i+1}$, $i = 1,2,3$ and $\sigma(C_4) = C_1$, and assume that $C_1$ is invariant by $\sigma^{4}$, then 
$\sigma^4(C_2) = \sigma^{4}(\sigma(C_1)) = \sigma(\sigma^4(C_1)) = \sigma(C_1) = C_2$. In particular if $C_1$ is pointwise  fixed, then one proves in a similar way that $C_2$ is pointwise fixed. A similar proof holds also for $C_3$ and $C_4$, so we get the statement.
\endproof

We denote by $2a_{\sigma}$ the number of smooth rational curves exchanged by $\sigma$ and  fixed by $\sigma^2$, and by $4s_{\sigma}$ the number of smooth rational curves cyclic permuted by $\sigma$ and pointwise fixed by $\sigma^4$ (and clearly they are interchanged by $\sigma^2$ two by two). Let now $a_{\sigma^{2}}$ be the number of the pairs of rational curves interchanged by $\sigma^2$ and pointwisely fixed by $\sigma^{4}$, then observe that $2a_{\sigma^{2}}=4s_{\sigma}$ and so $a_{\sigma^{2}} \in 2\mathbb{Z}$.


An important remark about the local behaviour of $\sigma^{2}$ in a neighborhood of a fixed point is the following:

\begin{remark}\label{rmk 3}
	The isolated fixed points by a non--symplectic automorphism $\sigma$ of type $(2,7)$ and $(3,6)$ are also isolated fixed points in $\Fix(\sigma^{2})$. The points of type $(4,5)$ in $\Fix(\sigma)$ are contained in a smooth fixed curve by $\sigma^{2}$. In fact the action of $\sigma^{2}$ at such points is given by the matrix $ \begin{pmatrix} 1 & 0 \\ 0 & \zeta_8^2 \end{pmatrix} $ which implies that these points belong to a smooth curve in $\Fix(\sigma^{2})$. For this reason we can say that if there exist points of $(4,5)$--type then $k_{\sigma^{2}} > k_{\sigma}$.
\end{remark} 

The following is a result about local actions of $\sigma$ at a point on a rational curve. With the same notation of Remark \ref{rmk 1} if $R$ is $\sigma$--invariant (not pointwise fixed) each action of $\sigma$ on $R$ has two fixed points. There are some restrictions about the possible actions of $\sigma$, in particular if we have an action on one of the fixed points then the action on the other point is determined. We recall the following result which is stated in a similar way in \cite[Lemma 8.1]{dillies2009order}. 
\begin{pro}\label{action}
	Let $\sigma$ be a non--symplectic automorphism of order $8$ on a K3 surface $X$ and suppose that the fixed locus of the involution $\sigma^{4}$ is the union of smooth rational curves. Then if $p_1$ is an isolated fixed point for $\sigma$ and it is contained on a rational curve $R$ fixed by $\sigma^4$, there exists another fixed point $p_2$ for $\sigma$ on $R$. If the local action in $p_1$ is of $(2,7)$--type then the local action in $p_2$ is of $(3,6)$--type and vice--versa. If the action in $p_1$ is of $(4,5)$--type then the action in $p_2$ is of $(4,5)$--type.
\end{pro}
\proof
We know that the action of a non--trivial finite order automorphism on a rational curve has two fixed points. Let $p_1$ and $p_2$ be the two fixed points on $R$. We can assume that the two fixed points are $p_1=(1:0)$ and $p_2=(0:1)$. Hence if $z$ is a local coordinate for $p_1$ on the rational curve $R\simeq \IP^1$ then $\frac{1}{z}$ is a local coordinate for $p_2$ on $R$. So that if $\sigma$ acts by $z\mapsto \zeta_8^j z$, $0\leq j\leq 7$ at $p_1$ then at $p_2$ 
it acts as  $\frac{1}{z}\mapsto \frac{1}{\zeta_8^j z}=\zeta_8^{8-j}\frac{1}{z}$. Observe that $j=0$ is not possible otherwise $R$ would be pointwise fixed by $\sigma$. Also $j=1,3,5,7$ are not possible. In fact $R$ is fixed by $\sigma^4$ so that we must have $\zeta_8^{4j}=(-1)^{j}=1$ which implies $j=2,4,6$. For these $j$'s one can compute that the local action at $p_1$ and $p_2$ (as points on $X$) can be diagonalized with $(t,s)$ as in the statement (up to permutation of the coordinates).
\endproof


\section{The classification}\label{classification_s}
Our goal is to give a classification of non--symplectic automorphisms of order eight under the assumption that their fourth power 
fixes only rational curves. 
Let $\sigma$ be such an automorphism, then
\begin{center}
	$\Fix(\sigma^4)=R'_1 \cup \dots \cup R'_T $
\end{center}
where the $R'_i$'s are smooth rational disjoint curves.
This implies that:

\begin{center}
	$\Fix(\sigma)=R_1 \cup \dots \cup R_{k_\sigma}\cup \{p_1, \dots, p_N\} $
\end{center}
which means that in the description of equation \eqref{luogofisso} we have $g(C)=0$.
 \\

\begin{theorem}\label{main k  1}
	Let $\sigma$ be a non--symplectic automorphism of order eight acting on a K3 surface $X$ with $\Pic(X)=S(\sigma^{4})$. Suppose that $\Fix(\sigma^{4})$ is not empty and it is the union of smooth rational curves. Then $k_{\sigma} \in \{0, 1\}$ and the invariants of $\sigma$ are as in Table \ref{tabella thm}. 
\end{theorem}
\proof 
Consider $p_1$ a fixed isolated point for $\sigma$. From Remark $\ref{rmk 1}$ there exists a smooth rational curve in $\Fix(\sigma^{4})$ such that $p_1 \in R'_i$.
	From Remark $\ref{rmk 2}$, $R'_i$ is $\sigma$--invariant. Since a finite order automorphism of a rational curve has two fixed points, there exists another fixed point for $\sigma$ on $R'_i$ : we call it $p_2$. We conclude that $N_{\sigma}$ is even. From what we know about the local behaviour of $\sigma$ at a fixed point, see Proposition \ref{action}, we know that if a fixed point is of $(2,7)$--type then the other fixed point on the same rational curve $R'_i$ is of $(3,6)$--type. If there is an action of $(4,5)$--type on $p_1$ then there is an action of $(4,5)$--type on $p_2$. By these considerations we obtain:
	$$n_{2,7}=n_{3,6}$$
	$$n_{4,5} \in 2\mathbb{Z}.$$
	Using Proposition $\ref{order 8}$ we obtain:
	$$\alpha=k_{\sigma}$$
	$$n_{3,6}=n_{2,7}=1+2k_{\sigma}\geq 1$$
	$$n_{4,5}=2+2k_{\sigma}\geq 2$$
Thanks to Remark \ref{rmk 3} the number of curves in $\Fix(\sigma^{2})$ is:
	$$ k_{\sigma^{2}}=k_{\sigma} + \frac{n_{4,5}}{2} +2a_{\sigma}.$$	
To conclude the proof we consider two cases, $k_{\sigma} \geq 1$ and $k_{\sigma}=0$.
	\begin{itemize}
	\item [$k_{\sigma} \geq 1$:]  
	In this case $k_{\sigma^{2}} \geq 3$ so the possible cases are $k_{\sigma^{2}}=3$ and $k_{\sigma^2}=4$ from Table \ref{g=0}. But from Remark \ref{rmk 0} we know that $m_{\sigma^{2}}$ has to be even, so checking again in Table \ref{g=0}, the unique case is $k_{\sigma^{2}}=3$ which has  $m_{\sigma^{2}}=2$. By using the previous equalities we get then that $n_{4,5}=4$, $k_{\sigma}=1$, $a_{\sigma}=0$ and $n_{4,5}=4$, 
$n_{2,7}=n_{3,6}=3$ so that $N_{\sigma}=n_{4,5}+n_{3,6}+n_{2,7}=10$. We can use Remark \ref{rmk 0} and Proposition \ref{order 8} to conclude that $r_{\sigma}=13$, $l_{\sigma}=3$, $m_{\sigma}=1$, $m_{1}=1$.\\ 
	\item [$k_{\sigma} =0$:] By using the equalities stated at the beginning of the proof we obtain:
	$$k_{\sigma} = 0,$$
	$$n_{3,6}=n_{7,2} =1,$$
	$$n_{4,5} = 2.$$
	So that in this case $N_{\sigma}=4$.
	Observe that $\sigma^{2}$ is an automorphism of order $4$ and it contains a rational curve in the fixed locus since two points of $(4,5)$--type in $\Fix(\sigma)$ are contained on a fixed curve for $\sigma^{2}$ so the number of curves in $\Fix(\sigma^{2})$ is:
	$$ k_{\sigma^{2}}=k_{\sigma} + \frac{n_{4,5}}{2} +2a_{\sigma}\geq 1.$$
	From Remark \ref{rmk 0} we know that $m_{\sigma^{2}}$ has to be even, so using Table $\ref{g=0}$ we conclude that there are four possible cases. \\
	If $\boldsymbol{m_{\sigma^{2}}=4}$ then $\boldsymbol{k_{\sigma}^{2}=1}$. From the previous equation $a_{\sigma}=0$. So we have  $(r_{\sigma^2},m_{\sigma^2},l_{\sigma^2},N_{\sigma^2},k_{\sigma^2},a_{\sigma^2})=(10,4,4,6,1,0)$.\\
	If $\boldsymbol{m_{\sigma^{2}}=2}$ then $k_{\sigma}^{2}\in\{3,2,1\}$. 
	If $\boldsymbol{k_{\sigma^{2}}=3}$ from the previous equation $a_{\sigma}=1$, so that $(r_{\sigma^2},m_{\sigma^2},l_{\sigma^2},N_{\sigma^2},k_{\sigma^2},a_{\sigma^2})=(16,2,2,10,3,0)$.\\
	If $\boldsymbol{k_{\sigma^{2}}=2}$ from the previous equation $a_{\sigma}=\frac{1}{2}$, which is not possible.\\
	If $\boldsymbol{k_{\sigma^{2}}=1}$ from the previous equation $a_{\sigma}=0$. In this case for $\sigma^{2}$ the invariants are $(r_{\sigma^2},m_{\sigma^2},l_{\sigma^2},N_{\sigma^2},k_{\sigma^2},a_{\sigma^2})=(12,2,6,6,1,2)$.\\
	We use then Remark \ref{rmk 0} and Proposition \ref{order 8} to compute the remaining invariants for $\sigma$ given in the Table \ref{tabella thm}.
\end{itemize}

\begin{table}[ht]
	\centering
\begin{tabular}  { |c|c|c|c|c|c|c|c|} 
	\hline
	$m_{1}$ & $m_{\sigma}$ & $r_{\sigma}$ & $l_{\sigma}$ & $N_{\sigma} $ & $k_{\sigma}$  & $a_{\sigma}$ & Examples \\
	\hline
	1 & 1 & 13 & 3 & 10 & 1 & 0 & \ref{k1}\\
	\hline
	2 & 2 & 6 & 4 & 4 & 0 & 0 & \\
	\hline
	1 & 1 & 9 & 7 & 4 & 0 & 1 &  \ref{k0}\\
	\hline
	1 & 3 & 7 & 5 & 4 & 0 & 0 &  \\
	\hline
	
\end{tabular} 
\vspace{0.2cm}
	
	\caption{Invariants of the automorphism}\label{tabella thm}
\end{table}

\endproof

\section{Elliptic fibrations}\label{elliptic_s}

\begin{definition}
	Let $X$ be a complex surface. An $\textbf{elliptic fibration}$ is a holomorphic map $f:X \longrightarrow B$ to a smooth curve $B$ such that the generic fiber is a smooth connected curve of genus one. A $\textbf{jacobian elliptic fibration}$ is an elliptic fibration admitting a section $ \pi: B \longrightarrow X$ such that $f \circ \pi=Id_{B}$. The surface $X$ is called an $\textbf{elliptic surface}$ if it admits an elliptic fibration
(not necessarily jacobian). We call $F_{v}$ the fiber $f^{-1}(v)$ over a point $ v \in B$.
	The $\textbf{Mordell--Weil group}$ is the group of sections of the elliptic fibration. 
\end{definition}
The $\textbf{zero section}$ of an elliptic fibration is a chosen section $s : B \longrightarrow X$ and we identify the map $s$ with the curve $s(B)$ on $X$. The point of intersection between the zero section and a fiber is the zero of the group law on the fiber. \\
For $K3$ surfaces we have that $B=\mathbb{P}^{1}$ (see \cite{miranda1989basic}) and, if the fibration is jacobian, it admits a Weierstrass equation:
\begin{equation}\label{eq W}
y^{2}=x^{3} + A(t)x + B(t),
\end{equation}
where $A(t)$ and $B(t)$ are two polynomials with complex coefficients with $t \in \mathbb{P}^{1}$  such that $\deg(A(t))=8$ and $\deg(B(t))=12$. Here the zero section is $ t \mapsto (0:1:0)$.\\
The discriminant of the fibration is a degree $24$ polynomial:
\begin{equation}\label{delta}
\Delta(t)=4A(t)^{3}+27B(t)^{2}.
\end{equation}
The equation \eqref{eq W} defines an elliptic fibration if and only if $\Delta(t)$ does not vanish identically.
Each zero of $\Delta(t)$ corresponds to a point $v$ of the base $\mathbb{P}^{1}$ such that $F_{v}$ is a singular fiber of the fibration. There are at most finitely many singular fibers. Let $\delta$ be the order of vanishing of $\Delta$ at a point corresponding to a singular fiber, by the $\textit{Kodaira classification}$ the possible singular fibers are recalled in Table $\ref{Kodaira class}$
where we denoted by $\Theta_{0}$ the component of a fiber meeting the zero section. The first column in the Table \ref{Kodaira class} contains the name of the reducible fiber according to Kodaira classification, the second the Dynkin diagram associated to the fiber, the last column contains the order of vanishing of $\Delta$ at the point corresponding to the singular fiber. \\
\begin{table}[ht]
\begin{eqnarray*}
\begin{array}{|c|c|c|c|}
\hline
\mbox{Name}&\mbox{Dynkin diagram}&\mbox{Description}&\delta\\
\hline
II&&\mbox{a cuspidal rational curve}& 2\\
\hline
III&\tilde A_1&\mbox{two rational curves meeting in a point of order 2}&3\\
\hline
IV&\tilde A_2&\mbox{three rational curves meeting at one point }&4\\
\hline
I_1&&\mbox{a nodal rational curve}&1\\
\hline
I_2&\tilde A_1&\tiny{\begin{array}{ccc}\Theta_0&=&\Theta_1\end{array}}&2\\
\hline
I_n&\tilde A_{n-1}&\tiny{\begin{array}{cccccc}\Theta_0&-&\Theta_1&-&\ldots&\Theta_i \\
|&&&&&|\\
\Theta_{n-1}&-&\Theta_{n-2}&-&\ldots&\Theta_{i+1}\end{array}}&n\\
\hline
I_k^*&\tilde D_{k+4}&\tiny{\begin{array}{cccccccccccc}
\Theta_0&&&&&&&&\Theta_{k+3}\\
&\diagdown&&&&&&\diagup&\\
&&\Theta_2&\ldots&\Theta_i -\Theta_{i+1} &\ldots&\Theta_{k+2}&\\
&\diagup&&&&&&\diagdown&\\
\Theta_1&&&&&&&&\Theta_{k+4}
\end{array}}&k+6\\
\hline
IV^*&\tilde E_6&\tiny{\begin{array}{ccccccccccccccccc}
\Theta_{0}&-&\Theta_{1}&-&\Theta_{2}&-&\Theta_{3}&-&\Theta_{4}\\
    & &     & &\!\mid & &  & &  & &   & \\
     & &      & & \Theta_5   & &   & &   & &   \\
     & &      & &\mid & &   & &   & &   &\\
     & &      & & \Theta_6   & &   & &   & &   \\
\end{array}}& 8\\
\hline
III^*&\tilde E_7&\tiny{\begin{array}{cccccccccccccc}
\Theta_{0}&-&\Theta_{2}&-&\Theta_{3}&-&\Theta_{4}&-&\Theta_{5}&-&\Theta_{6}&-&\Theta_{7}\\
      & &     & &&&\mid & &   & &   & &   & \\
      & &      & &&& \Theta_1   & &   & &   & &   & \\
\end{array}}& 9\\
\hline
II^*&\tilde E_8&\tiny{\begin{array}{ccccccccccccccc}
\Theta_0&-&\Theta_1&-&\Theta_2&-&\Theta_3&-&\Theta_4&-&\Theta_5&-&\Theta_6&-&\Theta_7\\
      & &      & &\mid& &   & &   & &    & &&&\\
      & &      & & \Theta_8   & &   & &   & &   & &&&
\end{array}}& 10\\
\hline
\end{array}
\end{eqnarray*}
\vspace{0.2cm}
	
	\caption{Kodaira classification}\label{Kodaira class}
\end{table}


A {\it simple component} of a fiber is a component with multiplicity one. In Table \ref{Dynkin_diagrams} we describe the singular fibers of an elliptic fibration with the multiplicities of the vertices of the extended Dynkin diagrams and we list the components with their multiplicities. The N\'{e}ron--Severi group of a surface admitting a jacobian elliptic fibration contains the class of a fiber $F$ and the class of the zero section $s$. Since the fibers are all algebraic equivalent, we have $ F \cdot F =0$. The zero section intersects each fiber in one point, so that $F \cdot s=1$.
The sections of an elliptic fibration on a $K3$ surface are smooth rational curves and this implies that their self--intersection is $-2$. Moreover, if $X$ is a $K3$ surface that admits a jacobian elliptic fibration, then there is an embedding of $\overline{U}$ in $\NS(X)$, where $\overline{U}$ is the two dimensional lattice

$$
\overline{U}=\left\lbrace \mathbb{Z}^{2},
\begin{pmatrix}
0 & 1 \\
1    & -2
\end{pmatrix} \right\rbrace.
$$
\\
Observe that by taking the generators $F$ and $F+s$ the lattice $\overline{U}$ is isometric to the hyperbolic two dimensional lattice

$$
U=\left\lbrace \mathbb{Z}^{2},
\begin{pmatrix}
0 & 1 \\
1    & 0
\end{pmatrix} \right\rbrace.
$$

If $f :X \longrightarrow \mathbb{P}^{1}$ admits an $n$-torsion section $s_n$ of order $n$ in the Mordell--Weil group then it induces an automorphism of the same order on $X$. This acts as the identity on the base of the fibration and as a translation by the section on each fiber, hence it is a symplectic automorphism \cite[Chapter 15, Lemma 4.4]{huybrechts2016lectures})

\begin{table}[ht]
  \centering

\renewcommand*{\arraystretch}{4.0}
\begin{tabular}{|c|c|p{4cm}|}
\hline
Name&simple components&Dynkin diagram \\ \hline
$\tilde A_n$& $\Theta_i,\,\, i=0,\ldots,n-1$&

\begin{picture}(10,30)
	\put(0,20){\circle*{2}}
	\put(10,20){\circle*{2}}
	\put(30,20){\circle*{2}}	
	\put(0,5){\circle*{2}}
	\put(10,5){\circle*{2}}
	\put(30,5){\circle*{2}}

	\put(-1,23){1}
	\put(9,23){1}
	\put(29,23){1}
	\put(-1,-3){1}
	\put(9,-3){1}
	\put(29,-3){1}
	
	\put(0,20){\line(1,0){10}}
	\put(0,20){\line(0,-1){15}}
	\put(30,20){\line(0,-1){15}}
	\put(0,5){\line(1,0){10}}
	
	\multiput(10,20)(4,0){5}{\line(1,0){2}} 
	\multiput(10,5)(4,0){5}{\line(1,0){2}} 
	
\end{picture}\\ \hline

$\tilde D_{k+4}$&$\Theta_i,\,\, i=0,1,k+3, k+4$&
\begin{picture}(10,30)
	\put(0,20){\circle*{2}}	
	\put(0,0){\circle*{2}}
	\put(10,10){\circle*{2}}
	\put(20,10){\circle*{2}}
	\put(40,10){\circle*{2}}
	\put(50,10){\circle*{2}}
	\put(60,20){\circle*{2}}
	\put(60,0){\circle*{2}}

	\put(0,20){\line(1,-1){10}}
	\put(0,0){\line(1,1){10}}
	\put(10,10){\line(1,0){10}}
	\put(40,10){\line(1,0){10}}
	\put(50,10){\line(1,1){10}}
	\put(50,10){\line(1,-1){10}}
	
\multiput(20,10)(4,0){4}{\line(1,0){2}} 
	
	\put(-1,22){1}
	\put(-1,2){1}
	\put(10,2){2}
	\put(20,2){2}
	\put(40,2){2}
	\put(49,2){2}
	\put(59,22){1}	
	\put(59,2){1}

\end{picture}\\ \hline

$\tilde E_6$&$\Theta_i,\,\, i=0,4,6$&\begin{picture}(10,30)
	\put(0,20){\circle*{2}}
	\put(10,20){\circle*{2}}
	\put(20,20){\circle*{2}}
	\put(30,20){\circle*{2}}
	\put(40,20){\circle*{2}}
	\put(20,10){\circle*{2}}
	\put(20,0){\circle*{2}}
	
	\put(0,20){\line(1,0){10}}
	\put(10,20){\line(1,0){10}}
	\put(20,20){\line(1,0){10}}
	\put(20,20){\line(0,-1){10}}
	\put(20,10){\line(0,-1){10}}
	\put(30,20){\line(1,0){10}}

	\put(-1,23){1}
	\put(9,23){2}
	\put(19,23){3}
	\put(29,23){2}
	\put(39,23){1}
	\put(23,10){2}	
	\put(23,0){1}

\end{picture}\\ \hline
$\tilde E_7$&$\Theta_i,\,\, i=0,7$&\begin{picture}(10,20)
	\put(0,10){\circle*{2}}
	\put(10,10){\circle*{2}}
	\put(20,10){\circle*{2}}
	\put(30,10){\circle*{2}}
	\put(40,10){\circle*{2}}
	\put(50,10){\circle*{2}}
	\put(60,10){\circle*{2}}
	
	\put(30,0){\circle*{2}}
	
	\put(0,10){\line(1,0){10}}
	\put(10,10){\line(1,0){10}}
	\put(20,10){\line(1,0){10}}
	\put(30,10){\line(0,-1){10}}
	\put(30,10){\line(1,0){10}}
	\put(40,10){\line(1,0){10}}
	\put(50,10){\line(1,0){10}}

	\put(-1,13){1}
	\put(9,13){2}
	\put(19,13){3}
	\put(29,13){4}
	\put(39,13){3}
	\put(49,13){2}
	\put(59,13){1}
	\put(33,0){2}	
\end{picture}\\\hline
$\tilde E_8$&$\Theta_7$ &\begin{picture}(10,20)
	\put(0,10){\circle*{2}}
	\put(10,10){\circle*{2}}
	\put(20,10){\circle*{2}}
	\put(30,10){\circle*{2}}
	\put(40,10){\circle*{2}}
	\put(50,10){\circle*{2}}
	\put(60,10){\circle*{2}}
	\put(70,10){\circle*{2}}

	\put(20,0){\circle*{2}}
	
	\put(0,10){\line(1,0){10}}
	\put(10,10){\line(1,0){10}}
	\put(20,10){\line(1,0){10}}
	\put(20,10){\line(0,-1){10}}
	\put(30,10){\line(1,0){10}}
	\put(40,10){\line(1,0){10}}
	\put(50,10){\line(1,0){10}}
	\put(60,10){\line(1,0){10}}

	\put(-1,13){2}
	\put(9,13){4}
	\put(19,13){6}
	\put(29,13){5}
	\put(39,13){4}
	\put(49,13){3}
	\put(59,13){2}
	\put(69,13){1}
	
	\put(23,0){3}

\end{picture}
\\\hline

\end{tabular}
\vspace{0.2cm}

\caption{Dynkin diagrams with the multiplicities of the components}\label{Dynkin_diagrams}
\end{table}


\section{Examples}\label{example_s}

We give here examples corresponding to the cases discussed in Theorem \ref{main k  1}. These are constructed by using jacobian elliptic fibrations on $K3$ surfaces. 
\subsection{Example $\boldsymbol{k_{\sigma}=1}$}\label{k1} 
The case $k_{\sigma}=1$ in Theorem \ref{main k  1} occurs, this means that we can find an example of a non--symplectic automorphism of order 8 on a K3 surface $X$ such that its fixed locus consists of a smooth rational curve and of 10 isolated points and the fixed locus of $\sigma^{4}$ consists of eight smooth rational disjoint curves. Consider the elliptic fibration on $X$ given by:
$$ y^{2}=x(x^{2}+tp_{6}(t))$$ with $p_6(t):=(a_6t^6+a_4t^4+a_2t^2+a_0)=(t^2-\alpha_1)(t^2-\alpha_2)(t^2-\alpha_3)$, where $a_6,a_4,a_2,a_0,\alpha_1,\alpha_2,\alpha_3 \in \mathbb{C}$, and the order 8 automorphism acting on it:
$$ \sigma:(x,y,t) \mapsto(-ix,\zeta_8y,-t) .$$

By \cite[Section 3]{kondo1992automorphisms} a holomorphic 2--form can be written as $\omega_X=\frac{dt \wedge dx}{2y}$ so that one computes
$$
\sigma(\omega_X)=\zeta_8 \omega_X, 
$$
hence $\sigma$ acts purely non--symplectically. Moreover $\sigma$ acts as an involution on the base $\mathbb{P}^{1}$ and it has order four on each fiber of the fibration.

The discriminant is $\Delta(t)=4t^3(t^2-\alpha_1)^3(t^2-\alpha_2)^3(t^2-\alpha_3)^3$.
Recall that $t \in\mathbb{P}^1$ so if we consider the homogenization of the polynomial in coordinates $(t:u)$, we obtain $\Delta(t,u)=4t^3(t^2-\alpha_1 u^2)^3(t^2-\alpha_2 u^2)^3(t^2-\alpha_3 u^2)^3u^3$.\\
We take now $\alpha_1=0$ and $\alpha_2=\alpha_3$. Under this assumtion the equation of the elliptic fibration becomes 
$$
y^{2}=x(x^{2}+t(t^2-\alpha_2)^{2})
$$ 
and $\Delta(t)=4t^9(t^2-\alpha_2)^6$ which in homogeneous coordinates $(t:u)$  is equal to $\Delta(t,u)=4t^9(t^2-\alpha_2u^2)^6u^3$. For generic choice of the coefficient $\alpha_2$ the fibration has $4$ singular fibers which correspond to the four zeros of $\Delta(t,u)$. To be more precise the fibration has a fiber of type $III^{*}$ over 0, which corresponds to the  $\tilde{E_7}$ Dynkin diagram, a fiber of type $III$ over $\infty$ which consists in two rational curves tangent in a point and two fibers of type $I_{0}^{*}$ over $\pm \sqrt{\alpha_2}$, see Table \ref{Kodaira class}. The action of $\sigma$ on the base fixes two points: $0$ and $\infty$, and so it preserves the fibers over these two points, the fibers $III$ and $III^{*}$, and it exchanges the two fibers of type $I_{0}^{*}$. If $f:X \longrightarrow \mathbb{P}^{1}$ is the elliptic fibration then $\Fix(\sigma)\subset f^{-1}(0) \cup f^{-1}(\infty)$.
Observe that the fibration has a two torsion section $s_{2}$, given by $t \mapsto (0:0:1)=(x:y:z)$ and the zero section $s$, given by $t \mapsto (0:1:0)=(x:y:z)$. These sections are preserved by the action of $\sigma$ and they have two fixed points corresponding to the intersection with the fibers over $0$ and $\infty$. These two sections are pointwise fixed by the action of $\sigma^2$. The sections $s_2$ and $s$ are simple sections which means that they meet one of the components of $III^*$ of multiplicity one (not the same), they meet one of the two components of $III$ (not the same) in a non--singular point and they meet the fibers of type $I_0^*$ in one of the components of multiplicity $1$ (not the same). We know that the components are not the same since  the sections $s$ and $s_2$ are contained in the fixed locus of $\sigma^4$ and the above mentioned components of the singular fibers are not fixed pointwise by $\sigma^4$ and cannot contain more than two fixed points. 
Now we can see that the two fixed points on each section are contained in $III^*$ and $III$. The fiber of type $III$ consists in two tangent rational curves. Each of these two rational curves has a fixed point which is not the tangency point, but a non--trivial finite order automorphism acting on a rational curve  has two fixed points, so we conclude that the double point on $III$ is fixed by $\sigma$. Observe that the sections $s$ and $s_2$ are preserved so that all the components of the fiber $III$ are preserved. On the fiber of type $III^*$ we have the two fixed points given by the intersections with the two sections $s$ ans $s_2$. Since the sections $s$ and $s_2$ are not exchanged then all components in $III^*$ are preserved. The unique component of this fiber of multiplicity four is preserved by $\sigma$ and it contains three fixed points so it is pointwise fixed.
 The component of multiplicity two which intersects the component of multiplicity four contains then a fixed point and we know that there is another fixed point on it. In conclusion $N_{\sigma}=10$ and $k_{\sigma}=1$. 

The square of the automorphism $\sigma^{2}:(x,y,t) \mapsto (-x,iy,t)$ preserves each fiber and acts as an automorphism of order four on it. Moreover $\sigma^{2}$ fixes two points on the generic smooth fiber, these two fixed points are contained in the two sections $s$ and $s_2$. This gives that $k_{\sigma^{2}} \geq 2$.\\
Since $l_{\sigma^{2}}=2m_{\sigma}=2$ and $m_{\sigma^{2}}=2m_1=2$ using Table $\ref{g=0}$ for the classification of non--symplectic automorphisms of order 4, we know that the invariants for $\sigma^{2}$ are $(r_{\sigma^{2}}, m_{\sigma^{2}},l_{\sigma^{2}}, N_{\sigma^{2}}, k_{\sigma^{2}}, a_{\sigma_{2}})=(16,2,2,10,3,0)$. The curves fixed by $\sigma^{2}$ are the curves fixed by $\sigma$ and the two sections $s$ and $s_{2}$.
The points fixed by $\sigma^{2}$ are 10 but they are not the same fixed points by $\sigma$, in fact the 4 fixed points for $\sigma$ on $s$ and $s_2$ now lie on fixed curves (namely the curves $s$ and $s_2$) and they do not give any contribution, but we add exactly 4 other points on the two fibers $I_0^*$. For this reason the number of fixed points remains the same.


Consider the curve defined by $y=0$. From the equation of the elliptic fibration we obtain $x=0$, which gives the zero section, and the curve $C$: $x^{2}+t^2(t^2-\alpha_2)^{2}=0$ which has a $2:1$ morphism to $\mathbb{P}^{1}$ and has possibly ramification points where $t^2(t^2-\alpha_2)^{2}=0$. These points lie on the four singular fibers: the fiber $III$, the two fibers of type $I_0^*$ and the fiber $III^*$. One has to study carefully the intersection of $C$ with the Dynkin diagram of the fibers to understand the ramification. The curve $C$ meets $III$ in the double point, $I_{0}^{*}$ in the two components of multiplicity one (so it does not ramify here) and $III^*$ in the component of multiplicity two.

 Recall that the Riemann--Hurwitz formula applied to a $2:1$ morphism from a curve $C$ to the projective line $\mathbb{P}^{1}$ is given by: 
 $$2g(C)-2=2(0-2)+\sum\limits_{p \in C } (e_p-1)$$
 where the sum runs over the ramification points which are two points in this case: the point on the fiber over 0, i.e. on $III^{*}$ and the point on the fiber over $\infty$ i.e on $III$ and $e_p$ is the ramification index at a ramification point $p$. In this case $e_p=2$. By using the formula 
 we can compute in an easy way the genus of the curve $C \subset \Fix(\sigma^{4})$ which is $g(C)=0$. Hence the fixed curves by $\sigma^{4}$ are all rational and we have eight of them : three components of the fiber $III^*$, two sections $s$ and $s_2$, the curve $C$ and two other rational curves which are fixed on the two fibers of type $I_0^*$.
%
 \subsection{Example $\boldsymbol{k_{\sigma}=0}$}\label{k0}
The case $k_{\sigma}=0$ in Theorem \ref{main k  1} occurs when $k_{\sigma^{2}}=3$. We can consider the same elliptic fibration of the previous example, and we fix $\alpha_1=0$ and $\alpha_2=\alpha_3$ as before. 
As we have already observed the fibration has a 2--torsion section given by $t \mapsto (0:0:1)=(x:y:z)$. Denote by $\tau$ the symplectic involution associated to this 2--torsion section. As we have observed before, this involution is symplectic. The involution exchanges the zero section $s$ and the 2--torsion section $s_2$. We cannot find fixed points for $\tau$ on the generic fiber since it acts as a translation, but we know (see \cite[Section 5]{Nikulin1}) that a symplectic involution has eight fixed points on a K3 surface. Consequently the 8 fixed points for $\tau$ are on the singular fibers of the elliptic fibration, recall that these are a fiber of type $III^*$, two fibers of type $I_0^*$ and a fiber of type $III$. 

On each of the two fibers of type $I_0^*$ we have two fixed points on the component of multiplicity two since the sections $s$ and $s_2$ are exchanged. On the fiber of type $III^*$ we have three fixed points, two of them are on the component of maximal multiplicity and the other is on the component of multiplicity two which intersects the component of maximal multiplicity. On the fiber of type $III$ the two rational curves are exchanged and so we have a fixed point for $\tau$ which is the double point. \\
We consider now $\sigma \circ \tau$. By a direct computation one computes $\sigma \circ \tau=\tau\circ\sigma$ so that we have an automorphism of order eight. Since $\tau$ is symplectic and $\sigma$ acts purely  non--symplectically we have that $\sigma \circ \tau$ acts purely non--symplectically too on the elliptic K3 surface. Now the fixed point for $\tau$ on $III$ is also a fixed point for $\sigma$, so it is a fixed point for $\sigma \circ \tau$. The fibers $I_0^*$ are exchanged by $\sigma$ so we cannot find fixed points for $\sigma \circ \tau$ on them and two of the three fixed points on $III^*$ are on a fixed curve for $\sigma$ so they contribute to $\Fix(\sigma \circ \tau)$.  Finally we conclude that $\Fix(\sigma \circ \tau)$ consists of four points. Where one of them is the double point of $III$ and three are on the fiber $III^*$. Finally since $\sigma$ and $\tau$ commute $(\sigma \circ \tau)^{2}=\sigma^2$, so that the behaviour of the order four automorphism is the same as it is described in the previous example.

\bibliographystyle{plain}
\bibliography{Biblio}

\begin{thebibliography}{10}

\bibitem{tabbaa2016order}
Dima Al~Tabbaa and Alessandra Sarti.
\newblock Order eight non-symplectic automorphisms on elliptic {$K3$} surfaces.
\newblock In {\em Phenomenological approach to algebraic geometry}, volume 116
  of {\em Banach Center Publ.}, pages 11--24. Polish Acad. Sci. Inst. Math.,
  Warsaw, 2018.

\bibitem{tabbaa2014classification}
Dima Al~Tabbaa, Alessandra Sarti, and Shingo Taki.
\newblock Classification of order sixteen non-symplectic automorphisms on {K}3
  surfaces.
\newblock {\em J. Korean Math. Soc.}, 53(6):1237--1260, 2016.

\bibitem{AS3}
Michela Artebani and Alessandra Sarti.
\newblock Non-symplectic automorphisms of order 3 on {$K3$} surfaces.
\newblock {\em Math. Ann.}, 342(4):903--921, 2008.

\bibitem{ASorder4}
Michela Artebani and Alessandra Sarti.
\newblock {Symmetries of order four on K3 surfaces}.
\newblock {\em J. Math. Soc. Japan}, 67(2):1--31, 2015.

\bibitem{ast}
Michela Artebani, Alessandra Sarti, and Shingo Taki.
\newblock {K3 surfaces with non-symplectic automorphisms of prime order.}
\newblock {\em Math. Z.}, 268(1-2):507--533, 2011.
\newblock with an appendix by Shigeyuki Kond{\=o}.

\bibitem{dillies2009order}
Jimmy Dillies.
\newblock On some order 6 non-symplectic automorphisms of elliptic {K}3
  surfaces.
\newblock {\em Albanian J. Math.}, 6(2):103--114, 2012.

\bibitem{dolgachev2007moduli}
Igor~V. Dolgachev and Shigeyuki Kond{\=o}.
\newblock Moduli of k3 surfaces and complex ball quotients.
\newblock In {\em Arithmetic and geometry around hypergeometric functions},
  pages 43--100. Springer, 2007.

\bibitem{huybrechts2016lectures}
Daniel Huybrechts.
\newblock {\em Lectures on K3 surfaces}, volume 158.
\newblock Cambridge University Press, 2016.

\bibitem{kondo1992automorphisms}
Shigeyuki Kondo.
\newblock Automorphisms of algebraic {K}3 surfaces which act trivially on
  {P}icard groups.
\newblock {\em Journal of the Mathematical Society of Japan}, 44(1):75--98,
  1992.

\bibitem{machida_oguiso}
Natsumi Machida and Keiji Oguiso.
\newblock On {$K3$} surfaces admitting finite non-symplectic group actions.
\newblock {\em J. Math. Sci. Univ. Tokyo}, 5(2):273--297, 1998.

\bibitem{miranda1989basic}
Rick Miranda.
\newblock {\em The basic theory of elliptic surfaces}.
\newblock ETS, 1989.

\bibitem{Nikulin1}
Viacheslav~V. Nikulin.
\newblock Finite groups of automorphisms of {K}\"ahlerian surfaces of type
  {$K3$}.
\newblock {\em Uspehi Mat. Nauk}, 31(2(188)):223--224, 1976.

\bibitem{nikulinintegral}
Viacheslav~V. Nikulin.
\newblock Integral symmetric bilinear forms and some of their applications.
\newblock {\em Math. USSR Izv.}, 14:103--167, 1980.

\bibitem{nikulinfactor}
Viacheslav~V. Nikulin.
\newblock Factor groups of groups of automorphisms of hyperbolic forms with
  respect to subgroups generated by 2-reflections. {A}lgebrogeometric
  applications.
\newblock {\em J. Soviet. Math.}, 22:1401--1475, 1983.

\bibitem{nikulindiscrete}
Viacheslav~V. Nikulin.
\newblock Discrete reflection groups in {L}obachevsky spaces and algebraic
  surfaces.
\newblock In {\em Proceedings of the {I}nternational {C}ongress of
  {M}athematicians, {V}ol. 1, 2 ({B}erkeley, {C}alif., 1986)}, pages 654--671.
  Amer. Math. Soc., Providence, RI, 1987.

\bibitem{rudakov1981surfaces}
Aleksei~Nikolaevich Rudakov and Igor~Rostislavovich Shafarevich.
\newblock Surfaces of type k3 over fields of finite characteristic.
\newblock {\em Itogi Nauki i Tekhniki. Seriya" Sovremennye Problemy Matematiki.
  Noveishie Dostizheniya"}, 18:115--207, 1981.

\bibitem{matthias}
Matthias Sch{\"u}tt.
\newblock {$K3$} surfaces with non-symplectic automorphisms of 2-power order.
\newblock {\em J. Algebra}, 323(1):206--223, 2010.

\bibitem{takiauto}
Shingo Taki.
\newblock Classification of non-symplectic automorphisms of order 3 on {K}3
  surfaces.
\newblock {\em Math. Nachr.}, 284:124--135, 2011.

\bibitem{Taki}
Shingo Taki.
\newblock Classification of non-symplectic automorphisms on {$K3$} surfaces
  which act trivially on the {N}\'eron-{S}everi lattice.
\newblock {\em J. Algebra}, 358:16--26, 2012.

\end{thebibliography}
\end{document}